\documentclass[12pt]{amsart}
\usepackage{amsfonts,amsmath,amstext,amsthm,amssymb}
\numberwithin{equation}{section}

\usepackage{hyperref,amsrefs}

\newtheorem{theorem}{Theorem}
\newtheorem{lemma}{Lemma}

\def\dist{{\rm dist\,}}
\def\supp{{\rm supp\,}}
\def\ZR{\ensuremath{\mathbb R}}

\def\ZI{\ensuremath{\mathbb I}}

\def\md#1#2\emd{\ifx0#1
\begin{equation*} #2 \end{equation*}\fi  
\ifx1#1\begin{equation}#2\end{equation}\fi   
\ifx2#1\begin{align*}#2\end{align*}\fi   
\ifx3#1\begin{align}#2\end{align}\fi    
\ifx4#1\begin{gather*}#2\end{gather*}\fi  
\ifx5#1\begin{gather}#2\end{gather}\fi   
\ifx6#1\begin{multline*}#2\end{multline*}\fi  
\ifx7#1\begin{multline}#2\end{multline}\fi  
}
\newcommand {\e }[1]{(\ref{#1})}
\newcommand {\lem }[1]{Lemma \ref{#1}}

\begin{document}
\title{An Estimate of the Maximal Operators\\ Associated with
 Generalized Lacunary Sets}

\author
{Grigor A. Karagulyan}

\address{Institute of Mathematics \\
Armenian National Academy of Sciences\\
Marshal Baghramian ave. 24b,\\
Yerevan, 375019, ARMENIA\\
}

\email{karagul@instmath.sci.am}

\author{Michael T. Lacey}

 \address{School of Mathematics \\
 Georgia Institute of Technology\\
 Atlanta GA 30332}
 
\email{lacey@math.gatech.edu}

\date{}
\maketitle

\begin{abstract}
Let $\Omega $ be any set of directions (unit vectors) on the plane. Denote by
$\mathcal{R}_\Omega $ the set of all rectangles which have a side parallel to
some direction from $\Omega $. In this paper we study maximal operators on
the plane $\ZR^2$ defined by 
\md0
M_\Omega f(x)=\sup_{x\in R\in \mathcal{R}_\Omega}\frac{1}{|R|}\int_R|f(y)|dy.
\emd
We are interested in extensions of lacunary sets of directions, to collections 
we call $N$--lacunary, for integers $N$. We proceed by induction. Say that 
$\Omega=\{v_k\mid k\in\mathbb N \}$ is $1$--lacunary iff for each integer $k$, $v_k$ and $v_{k+1}$ are 
neighboring  points, and  there is a direction $v_\infty$ so that 
\md0
\tfrac12|v_k-v_{k+1}|<|v_{k+1}-v_\infty|<|v_k-v_{k+1}|.
\emd 
Every $N+1$--lacunary set can be obtained from some $N$--lacunary
$\Omega_N$ adding some points to $\Omega_N$. Between each two
neighbor points $a,b\in\Omega_N$ we can add a  $1$--lacunary
sequence (finite or infinite).    We show that for all $N$ lacunary sets 
$\Omega$, 
$$
\|M_\Omega f(x)\|_2\lesssim{}N \|f\|_2.
$$
Observe that every set $\Omega$ of $N$ points is $(C\log N)$--lacunary.  We 
then obtain a Theorem of N.~Katz \cite{Katz2}.
Both the current inequality, and Katz' result are  consequence of a general result 
of Alfonseca, Soria, and Vargas \cites{ASV2}.  We offer the current proof as a succinct, 
self--contained approach to this inequality. 

\end{abstract}

\begin{section}{Introduction}\everypar{} \parskip=11pt

Let $\Omega $ be any set of directions (unit vectors) on the plane. Denote by
$\mathcal{R}_\Omega $ the set of all rectangles which have a side parallel to
some direction from $\Omega $. In this paper we study maximal operators on
the plane $\ZR^2$ defined by 
\md1\label{MOm}
M_\Omega f(x)=\sup_{x\in R\in \mathcal{R}_\Omega}\frac{1}{|R|}\int_R|f(y)|dy.
\emd
A.~Nagel, E.M.~Stein and S.~Wainger \cite{NSW} using Fourier transform method
proved the boundedness of $ M_\Omega f(x)$ in spaces $L^p$, $1<p<\infty $
for any lacunary set of directions ${\Omega }=\{\theta_k\}$,
$(\arg\theta_{k+1}<\lambda\arg \theta_k, \,\, \lambda <1)$. 

We are interested in extensions of lacunary sets of directions, to collections 
we call $N$--lacunary, for integers $N$. We proceed by induction. Say that 
$\Omega=\{v_k\mid k\in\mathbb N \}$ is $1$--lacunary iff for each integer $k$, $v_k$ and $v_{k+1}$ are 
neighboring  points, and  there is a direction $v_\infty$ so that 
\md0
\tfrac12|v_k-v_{k+1}|<|v_{k+1}-v_\infty|<|v_k-v_{k+1}|.
\emd 
Every $N+1$--lacunary set can be obtained from some $N$--lacunary
$\Omega_N$ adding some points to $\Omega_N$. Between each two
neighbor points $a,b\in\Omega_N$ we can add a  $1$--lacunary
sequence (finite or infinite). So if $\Omega$ is some $N$--lacunary set we can fix
a sequence of sets $\Omega_1\subset\Omega_2\subset\cdots
\subset\Omega_{N-1}\subset\Omega $ such that each $\Omega_k$ is $k$--lacunary.


It is commonly known that maximal functions in $N$--lacunary directions are bounded for all integers $N$.
For instance, the case of $2$--lacunary is due to P.~Sj\"ogren and P.~Sj\"olin \cite{SS}. 
We are interested in growth of  the norm of $M_\Omega$ for $N$--lacunary, as $N$ tends to 
infinity. 

\begin{theorem}  \label{Tlac}
For all integers $N$, 
and all $N$--lacunary sets $\Omega$ we have
$$
\|M_\Omega f(x)\|_2\lesssim{}N \|f\|_2.
$$
\end{theorem}

It is easy to check that each set of directions of cardinality $N$ is $(C\log N)$--lacunary, for 
an absolute constant $C$.  Therefore, as a corollary, we see that for finite collections $\Omega$, we have 
\md1 \label{e.katz}
\| M_\Omega f\|_2\lesssim{}(\log \sharp\Omega)\|f \|_2.
\emd
This inequality is due to N.~Katz \cite{Katz2}.  This estimate is sharp as the power of $(\log \sharp\Omega)$, and so 
in the  Theorem, our estimate is sharp as to the power of $N$.  

Both Katz' result and our Theorem is a consequence of a more general result of Alfonseca, Soria, and Vargas \cite{ASV2}, a result we recall in more detail below.  The current proof is succinct, and self--contained, and so  may prove to be of some independent interest. 

\bigskip 

We close this section with a more detailed, but far from complete,
 description of the history of this question, and the relationship of our result to the literature.  
 In 1977, 
A.~Cordoba \cite{MR56:6259} considered the maximal function formed over all rectangles that are 
$1$ by $N$, obtaining a slow increase in the norm on $L^2$.  Thus, the set $\Omega$ is uniformly distributed, 
but one only considers rectangles of one aspect ratio.   The method of proof employed a  
geometric method to prove a covering lemma.  The method, as  described in A.~Cordoba and R.~Fefferman \cite{CoFe}, was broadly influential.  The point of 
view adopted in this paper was formalized in an article from 1979 by S.~Wainger \cite{MR82g:42018}. 
The estimate (\ref{e.katz}) in the instance of uniformly distributed directions was proved by 
J.~Stromberg \cite{Str2}, in 1978.  

On the other hand, there were natural reasons to expect that the instance of lacunary directions would 
behave differently, and was investigated by J.~Stromberg \cite{Str1}.  The full range of $L^p$, $1<p<\infty$, 
inequalities in this instance was established by Fourier analysis, and square function methods by 
A.~Nagel, S.~Wainger, and E.M.~Stein \cite{NSW}, a method that also proved to be influential. 
These results are related to interesting results on multipliers, as shown by  A.~Cordoba and R.~Fefferman 
\cite{MR55:6096}.  For extensions of this, see 
A.~Carbery \cite{MR89h:42026}.

An interesting question was if  Stromberg's result \cite{Str2} in the uniformly distributed case 
extended to the case of $N$ distinct directions.  A partial result was treated by Barrionuevo 
\cites{MR93f:42038,Barr}.  And the definitive result was obtained by N.~Katz \cite{Katz2}. His method of proof 
is a clever duality argument, relying on an John--Nirenberg type to obtain the required estimate.

At this point, we note that there is a distinction between the case of rectangles of all aspect 
ratios, as we do, and the case of a fixed aspect ratio.  It is the later case that is considered by 
e.g.~A.~Cordoba \cite{MR56:6259}, and in Katz' paper \cite{Katz1}.

An interesting question concerns the maximal function computed in a set of directions specified by 
a Cantor set of directions.  For the ordinary middle third Cantor set, there is a partial result on $L^2$ 
by A.~Vargas \cite{MR96m:42033}.  Yet, this full maximal function is unbounded on $L^2$, as proved  
by N.~Katz \cite{MR98b:42032}.  It would be interesting to obtain meaningful information about this maximal 
operator on $L^p$, for $p>2$.  K.~Hare \cite{MR2003f:42027} uses Katz' argument, with more general Cantor sets.

Recently, A.~Alfonesca, F.~Soria and A.~Vargas \cites{ASV, ASV2}, also see Alfonseca \cite{Alf},  have proved an interesting orthogonality principle 
for these maximal functions.  Let $\Omega=\{v_k\mid k\in\mathbb N\}$ be a set of directions, and between 
two neighboring directions $v_k,v_{k+1}$, let $\Omega_k$ be an arbitrary set of directions.  Then, (\cite{ASV2}) it is the case that 
\md0
\| M \|_{2\to2} \le{}C\| M_\Omega \|_{2\to2}+\sup_k \|M_{\Omega_k} \|_{2\to2}.
\emd
What is essential is that the second term occurs with constant $1$. 
This proves our Theorem.  Let $\eta(N)$ be the maximum of $\| M_{\Omega_N}\|_{2\to2}$, with the maximum taken over all $N$--lacunary sets of directions.  The inequality above clearly implies that 
$\eta(N)\le{}C\eta(1)+\eta(N-1)$.  Iterating the inequality $N-1$ times proves the Theorem.   

General necessary and sufficient conditions on $\Omega$ for the boundedness of $M_\Omega$ have been sought by 
J.~Duoandikoetxea, and A.~Vargas \cite{MR97c:42031}, with extensions by K.~Hare,  and J.~R{\"o}nning \cites{MR2000a:42035, MR2003b:42032}.

A paper by M.~Christ \cite{MR92f:42024} includes examples of sets of directions $\Omega$, and partial results on the norm boundedness of 
$M_\Omega$ which are not incorporated into the theories associated with this subject.   K.~Hare and F.~Ricci \cite{MR1938749} have 
established an interesting variant of the lacunary directional maximal function.

\end{section}


\begin{section}{Notations}  \everypar{} \parskip=11pt

By $A\lesssim{}B$ we mean that there is an absolute constant $K$ so that $A\le{}KB$.   By $\widehat f(\xi)$, we mean the 
Fourier transform of $f$, thus  
\md0
\widehat f(\xi )=\int f(x)e^{ix \cdot \xi }\; dx  
\emd

\medskip 

We use a well--known reduction to parallelograms.  It is clear that we can associate directions in $\Omega$ to points in 
e.g.~$(0,1/4)$.   Denote
\md1\label{Mal}
P_\alpha f(x)=\sup_{\delta_1,\delta_2}\frac{1}{4\delta_1\delta_2 }
\int_{x_1-\delta_1 }^{x_1+\delta_1 }\int_{x_2-x_1\alpha -\delta_2 }^{x_2-x_1\alpha +\delta_2 }
|f(t_1,t_2)|\; dt_2dt_1.
\emd
This is a maximal function over parallelograms, with one side parallel to the $x$ axis, and the other side forming 
an angle of slope $\alpha$ with the $x$ axis.  
Then in order to prove the theorem it is sufficient to prove
\md0
\|\sup_{\alpha \in\Omega }P_\alpha f\|_2\le CN\|f\|_2
\emd
where $\Omega $ is any $N$--lacunary set from $(0,1)$.

Our method of proof is Fourier analytic, and we shall find it convenient to use the 
 the Fejer kernel
\md0
K_r(x)=\int_{-r}^r\left(1-\frac{|t|}{r}\right)e^{-itx}dt=
\frac{}{}\frac{4\sin^2\frac{Nx}{2}}{Nx^2}
\emd
For any $r,R$ with $0\le r<R/2$ we define the following functions
\md0
\psi_r(x)=2K_{2r}(x)-K_{r}(x),\quad
\psi_{r,R}(x)=\psi_R(x)-\psi_r(x)
\emd
Sometimes we will write $\psi_{0,r}$ instead of $\psi_r(x)$. We have
\md1\label{lin}
\widehat\psi_{r,R}(\xi)={}
\left\{
\begin{array}{rcl}
1 &\hbox{ if }& |\xi|\in [2r,R]\\
0& \hbox{ if }& 0\le |\xi|\le r \hbox { or } |\xi|>2R\\
\hbox {linear } &\hbox{ on each }& \pm [r,2r],\pm [R,2R]
\end{array}
\right.
\emd
From a property of Fejer kernel we have
\md0
|\psi_{r,R}(x)|\le C\left(\max\left\{\frac{1}{Rx^2},R\right\}+
\max\left\{\frac{1}{rx^2},r\right\}\right)
\emd
Thus for some sequence of intervals $\omega_k=\omega_{k,r,R}$ with centers at  $0$.
\md5\label{int} \begin{split}
|\psi_{r,R}(x)|\le C\sum_k\gamma_k\frac{\ZI_{\omega_k}(x)}{|\omega_k|}=
 \zeta_{r,R}(x)\\
\gamma_k>0,\qquad \sum_k\gamma_k<1, \qquad \omega_k\supset(1/R,1/R) .
\end{split}
\emd

Choose a Schwartz function $\phi $ with
\md1
\phi \ge 0,\quad
\supp\widehat\phi \subset[-1,1]\label{psi}.
\emd
We can fix an even function $\lambda $ with
\md1
\max\{|\phi (x)|,|x\phi (x)|\}\le \lambda (x),\qquad
\int_\ZR\lambda (x)dx\le C,\label{xi}
\emd

Then define a Fourier analog of the average over parallelograms by 
\md1\label{Gam}
\Gamma_{r,R,h}^\alpha f(x)=\big(\psi_{r,R}(x_2-x_1\alpha)\phi_h(x_1)\big)*
f(x),\quad x=(x_1,x_2)\in \ZR^2.
\emd
where
\md0
\phi_h(x)=\frac{1}{h}\phi \left(\frac{x}{h}\right).
\emd
From \e{Gam} and \e{Mal} it follows that
\md0
P_\alpha f(x)\le C\sup_{R,h}\Gamma_{R,h}^\alpha f(x).
\emd
and therefore to prove our Theorem we need to verify the inequality 
\md1\label{the}
\|\sup_{R,h,\alpha\in\Omega}\Gamma_{R,h}^\alpha f(x)\|_2\le CN\|f\|_2
\emd

Taking the Fourier transform both sides of \e{Gam} we get
\md1\label{GamF}
\widehat\Gamma_{r,R,h}^\alpha f(\xi)=
\widehat\phi(h(\xi_1+\xi_2\alpha ) )\widehat\psi_{r,R}(\xi_2)\widehat f(\xi)
\emd

\end{section}


\begin{section}{Proof of Theorem}\everypar{} \parskip=11pt

\begin{lemma}\label{L1}
Let $\alpha,\beta\in (0,1)$ be any numbers and
$0<r<R,h>0$. The operator $\Gamma_{r,R,h}^\alpha f(x)$ defined in \e{Gam}
satisfies pointwise estimate
\md1\label{lem1}
|\Gamma_{r,R,h}^{\alpha}f(x)|
\le C\left(h R|\alpha-\beta|+1\right)P_{\beta}f(x),
\quad x\in \ZR^2.
\emd
\end{lemma}

\begin{proof} From \e{int} we have
\md0
\psi_{r,R}(x_2-x_1\alpha)\le C\sum_k\frac{\gamma_k}{|\omega_k|}\ZI_{\omega_k}
(x_2-x_1\alpha)
\emd
where we have $|\omega_k|>2/R$.
Denote $\lambda(x_1)=2Rx_1|\alpha-\beta|+2$ and assume
\md1\label{o1}
x_2-x_1\alpha\in \omega_k
\emd
for some $k$. Then taking account of \e{int} we get
\md5 \nonumber
\left|\frac{x_2-x_1\beta}{\lambda(x_1)}\right|=
\left|\frac{x_2-x_1\alpha+x_1(\alpha-\beta)}{\lambda(x_1)}\right|\\  \label{Rbig?}
\le\left|\frac{x_2-x_1\alpha}{2}\right|+
\frac{1}{2R}\le\frac{|\omega_k|}{2},
\emd
which means
\md1\label{o2}
\frac{x_2-x_1\beta}{\lambda(x_1)}\in \omega_k.
\emd
Hence we conclude that \e{o1} implies \e{o2}. Therefore
\md0
\ZI_{\omega_k}(x_2-x_1\alpha)\le
\ZI_{\omega_k}\left(\frac{x_2-x_1\beta}{\lambda(x_1)}\right)
\emd
Finally we get
\md4
\psi_{r,R}(x_2-x_1\alpha)\le C\sum_k\frac{\gamma_k}{|\omega_k|}\ZI_{\omega_k}
\left(\frac{x_2-x_1\beta}{\lambda(x_1)}\right)\le
\zeta_{r,R}\left(\frac{x_2-x_1\beta}{\lambda(x_1)}\right)
\emd
Thus taking account of \e{xi} we obtain
\md2
\frac{1}{h}\phi \left(\frac{x_1}{h}\right)
\zeta_{r,R}\left(\frac{x_2-x_1\beta}{\lambda(x_1)}\right)
{}\le{}
 C\left(h R|\alpha-\beta|+1\right)\frac{1}{h}\xi\left(\frac{x_1}{h}\right)
\frac{1}{\lambda(x_1)}\zeta_{r,R}(\frac{x_2-x_1\beta}{\lambda(x_1)})
\emd
from which we easily get \e{lem1}. \end{proof}

For any interval $J=(a,b)$ we denote by $S(J)$ the sector
$\{ax_2\le x_1\le bx_2\}$. For any sector $S$ define by $2S$ the sector
which has same bisectrix  with $S$ and twice bigger angle. Denote
by $T_Sf$ the multiplier operator defined
$\widehat T_Sf=\ZI_S\widehat f$.


\begin{lemma}\label{L2}
Let $J_1\supset J_2\supset\cdots \supset J_n$
be some sequence of intervals with
\md1
J_k=[\alpha_k,\beta_k]\subset (0,1),\quad
\dist((J_k)^c,J_{k+1})\le |J_{k+1}|, \quad 1\le{}k\le{}n\label{J}
\emd
Then for any $\theta \in \bigcap J_k$ and any function $f\in L^2(\ZR^2)$
we have
\md3 \begin{split} \label{lem2}
P_\theta f\lesssim{}& P_0f+P_{\theta }(T_{2S(J_n)}f)\\
{}&\quad{}+\sum_{k=1}^{n-1}
P_{\alpha_k}(T_{2S(J_k)}f)+P_{\beta_k}(T_{2S(J_k)}f)
\end{split}
\emd
where $P_0$ is a $P_\alpha $ with $\alpha =0$.
\end{lemma}

\begin{proof}

  Regard $\theta\in \bigcap J_k$ as fixed. For any $R,h$ we have
\md1\label{GamR}
\widehat\Gamma_{R,h}^\theta f(\xi)=
\widehat\psi_R(\xi_2)\widehat\phi(h(\xi_1+\xi_2\theta))\widehat f(x)
\emd
Denote
\md1\label{rk}
r_0=0,\quad r_k=\frac{2}{h|J_k|}\quad 1\le{}k\le{}n.
\emd
From \e{lin} it follows that
\md1\label{psis}
\widehat\psi_R(\xi_2)=\sum_{k=1}^m\widehat\psi_{2r_{k-1},r_k}(\xi_2)+
\widehat\psi_{2r_m,R}(\xi_2)
\emd
 where $ m=\max\{ k:r_k<2R\}$.  Denote
\md2
\Gamma_kf(x)={}&\Gamma_{2r_k,r_{k+1},h}^\theta f(x)\quad 0\le{}k<m,
\\
 \Gamma_mf(x)={}&\Gamma_{2r_m,R,h}^\theta f(x).
\emd
Then by \e{GamF} we have
\md2
\widehat\Gamma_k f(\xi)
={}&\widehat\psi_{2r_{k-1},r_k}(\xi_2)\widehat\phi(h(\xi_1+\xi_2\theta))\widehat f(x)
\quad 1\le{}k<m\\
\widehat\Gamma_m f(x)
={}&\widehat\psi_{2r_m,R}(\xi_2)\widehat\phi(h(\xi_1+\xi_2\theta))\widehat f(x)
\emd
and therefore using \e{psis} we obtain
\md1\label{gsum}
\Gamma_{R,h}^\theta f=\sum_{k=0}^m\Gamma_k f
\emd

\everypar{} \parskip=11pt

Let us show
\md3 \begin{split}\label{supp}
\supp\widehat\psi_{2r_k,r_{k+1}}(\xi_2)\widehat\phi(h(\xi_1+\xi_2\theta))\subset{}&
2S(J_k),\quad 1\le{}k<m,\\
\supp\widehat\psi_{2r_m,R}(\xi_2)\widehat\phi(h(\xi_1+\xi_2\theta))\subset{}&2S(J_m)
\end{split}
\emd
From which it follows that 
\md0
\Gamma_k f=\Gamma_k \big(T_{2S(J_k)}f\big), \quad 1\le{}k\le{}m
\emd
Indeed, from \e{psi} and \e{lin} it follows that
\md4
\supp\widehat\psi_{2r_k,r_{k+1}}(\xi_2)\widehat\phi(h(\xi_1+\xi_2\theta))\\
=\{(\xi_1,\xi_2): r_k\le \xi_2\le 2r_{k+1},\,\, |\xi_1+\xi_2\theta |<\frac{1}{h}\}
\emd
The last set is a parallelogram with vertexes
$(r_k\theta\pm \frac{1}{h},r_k)$ and $(2r_{k+1}\theta\pm \frac{1}{h},2r_{k+1})$.
These vertexes are from $2S(J_k)$ because
\md0
\frac{r_k\theta\pm \frac{1}{h}}{r_k}=\theta\pm \frac{|J_k|}{2}
\emd
which means  $(r_k\theta\pm \frac{1}{h},r_k)\in 2S(J_k)$. The same conclusion
is true for next the pair of vertexes. This implies \e{supp}.

Using \lem{L1} we conclude
\md3 \begin{split}\label{GT}
|\Gamma_k f|
{}\lesssim{}& (hr_{k+1}\min\{|\theta-\alpha_k|,|\theta-\beta_k|\}+1)\times
 \\&\quad(P_{\alpha_k}\big(T_{2S(J_k)}f\big)
 {}+{}
P_{\beta_k}\big(T_{2S(J_k)}f\big))\quad 1\le{}k<m
\end{split}
\emd
Notice also
\md5
|\Gamma_0 f|\le{} P_0f\label{M0}\\
|\Gamma_m f|\le{} P_\theta T_{2S(J_m)}f\label{Mm}
\emd
By $\theta\in J_{k+1}\subset J_k$ and \e{J} we have
\md0
\min\{|\theta-\alpha_k|,|\theta-\beta_k|\}\le 2|J_{k+1}|
\emd
The last with \e{rk} implies
\md0
hr_{k+1}\min\{|\theta-\alpha_k|,|\theta-\beta_k|\}\le 4
\emd
Hence by \e{GT} we observe
\md0
|\Gamma_k f|
{}\lesssim{} P_{\alpha_k}\big(T_{2S(J_k)}f\big)+
P_{\beta_k}\big(T_{2S(J_k)}f\big),\quad 1\le{}k<m.
\emd
Finally taking account also \e{M0} and \e{Mm} we get \lem{L2}.

\end{proof}


\begin{proof}[Proof of Theorem~\ref{Tlac}]

Let $\Omega \subset (0,1)$ be any N-lacunary set. We fix the sets
$\Omega_1\subset\Omega_2\subset\cdots\subset\Omega_{N-1}\subset\Omega_N=\Omega $
from definition of N-lacunarity.
Fix any angle $\theta \in \Omega $ and $R,h>0$. Suppose
\md1
\theta\in \Omega_m\setminus \Omega_{m-1},\hbox { for some } m\le N.
\emd
Denote by $G_k$ the set of all intervals whose vertexes are neighbor
points in $\Omega_k$. We can choose a sequence of intervals
$J_k=[\alpha_k,\beta_k]\in G_k$ $k=1,2,\cdots ,m$ such that
\md0
\theta\in \bigcap_{1\le k\le m} J_k, \qquad
\theta=\alpha_m \quad \text{(or  $\theta=\beta_m$)}
\emd
It is clear that sequence $J_k$ satisfies conditions of \lem{L2}.
Hence,
\md2
|M_\theta f|^2\lesssim{}& \big\{M_0f+\sum_{k=1}^m
(M_{\alpha_k}(T_{2S(J_k)}f)+M_{\beta_k}(T_{2S(J_k)}f))\big\}^2
\\
{}\lesssim{}&  |M_0f|^2+m\sum_{k=1}^m 
|M_{\alpha }(T_{2S(J )}f)|^2+|M_{\beta }(T_{2S(J )}f)|^2 
\emd
and therefore, summing over every interval $J=(\alpha,\beta)\in G_k$,
\md1  \label{log}
\sup_{\theta\in \Omega }|M_\theta f|^2 
{}\lesssim{}  |M_0f|^2+N\sum_{k=1}^N\sum_{J=(\alpha , \beta )\in G_k}
|M_{\alpha }(T_{2S(J )}f)|^2+|M_{\beta }(T_{2S(J )}f)|^2 
\emd
On the other hand using the $(2,2)$ bound of strong maximal operator we get for each $1\le{} k\le N$, 
\md2
\int_{\ZR^2} \sum_{J=(\alpha , \beta )\in G_k}
|M_{\alpha }(T_{2S(J )}f)|^2+|M_{\beta }(T_{2S(J )}f)|^2dx
{}\lesssim{} & \int_{\ZR^2} \sum_{J=(\alpha , \beta )\in G_k}\ZI_{2S(J)}|\widehat f|^2 d\xi
\\{}\lesssim{}& \int_{\ZR^2} |\widehat f|^2d\xi
\\{}={}&\int_{\ZR^2} |f|^2dx
\emd
Finally taking account of \e{log} we obtain
\md0
\int_{\ZR^2} \sup_{\theta\in \Omega }|M_\theta f|^2dx\lesssim{}
N^2 \int_{\ZR^2} |f|^2dx
\emd
\end{proof}

\end{section}

\end{document}